\documentclass[preprint,12pt]{elsarticle}

\usepackage{amsmath,amsfonts,latexsym,subfigure,amsthm,amssymb,wasysym}
\usepackage{graphicx}
\usepackage{tikz}
\usepackage{ifthen}
\usepackage[usenames,dvipsnames]{pstricks}
\usepackage{epsfig}
\usepackage{pst-grad} 
\usepackage{pst-plot} 
\usepackage{enumitem}
\usepackage{algorithm}
\usepackage{algorithmic}
\usepackage{wrapfig}
\usepackage{color}

\def\calN{{\mathcal N}}
\def\calD{{\mathcal D}}

\usepackage{bbm}

\DeclareMathOperator*{\argmax}{argmax}

\theoremstyle{remark}

\journal{Elsevier}

\begin{document}

\begin{frontmatter}



\title{Reduced Basis Decomposition: a Certified and Fast Lossy Data Compression Algorithm}


\author[UMassD]{Yanlai Chen \fnref{funding}}
\ead{{yanlai.chen@umassd.edu}}
\ead[url]{www.faculty.umassd.edu/yanlai.chen/}
\address[UMassD]{Department of Mathematics, University of Massachusetts Dartmouth, 285 Old Westport Road, North Dartmouth, MA 02747, USA.}

\fntext[funding]{This research was partially supported by National Science Foundation grant DMS-1216928.}

\begin{abstract}
Dimension reduction is often needed in the area of data mining. The goal of these methods is to map the given high-dimensional data into a low-dimensional space preserving certain properties of the initial data. There are two kinds of techniques for this purpose. The first, projective methods, builds an explicit linear projection from the high-dimensional space to the low-dimensional one. On the other hand, the nonlinear methods utilizes nonlinear and implicit mapping between the two spaces. In both cases, the methods considered in literature have usually relied on computationally very intensive matrix factorizations, frequently the Singular Value Decomposition (SVD). The computational burden of SVD quickly renders these dimension reduction methods infeasible thanks to the ever-increasing sizes of the practical datasets. 

In this paper, we present a new decomposition strategy, Reduced Basis Decomposition (RBD), which is inspired by the Reduced Basis Method (RBM). Given $X$ the high-dimensional data, the method approximates it by $Y \, T (\approx X)$ with $Y$ being the low-dimensional surrogate and $T$ the transformation matrix. $Y$ is obtained through a greedy algorithm thus extremely efficient. In fact, it is significantly faster than SVD with comparable accuracy. $T$ can be computed on the fly. Moreover, unlike many compression algorithms,  it easily finds the mapping for an arbitrary ``out-of-sample'' vector and it comes with an ``error indicator'' certifying the accuracy of the compression. Numerical results are shown validating these claims.

\end{abstract}

\begin{keyword}

Data mining, Lossy compression, Reduced basis method, Singular value decomposition, Greedy algorithm



\end{keyword}

\end{frontmatter}



\section{Introduction}
\label{sec:intro}

Dimension reduction is ubiquitous in many areas ranging from pattern recognition, clustering, classification, to fast 
numerical simulation of complicated physical phenomena. The fundamental question to address is how to approximate 
a $n$-dimensional space by a $d$-dimensional one with $d \ll n$. Specifically, we are given a set of high-dimensional data 
\begin{equation}
\label{eq:highd_x}
X = [x_1, x_2, \dots, x_n] \in {\mathbb R}^{m \times n},
\end{equation}
and the goal is to find its low-dimensional approximation 
\begin{equation}
\label{eq:lowd_y}
Y = [y_1, y_2, \dots, y_d] \in {\mathbb R}^{m \times d}
\end{equation}
with reasonable accuracy.

There are two types of dimension reduction methods. The first category consists of ``projective'' ones. These are the 
linear methods that are {\em global} in nature, and that explicitly transform the data matrix $X$ into a low-dimensional one by $Y = T X$. 
The leading examples are the Principal Component Analysis (PCA) and its variants. The methods in the second category act 
locally and are inherently nonlinear. For each sample in the high-dimensional space (e.g. each column of $X$), they directly find 
their low-dimensional approximations by preserving certain locality or affinity between nearby points.

In this paper, inspired by the reduced basis method (RBM), we propose a linear method called ``Reduced Basis Decomposition (RBD)''. 
It is much faster than PCA/SVD-based techniques. Moreover, its low-dimensional vectors are equipped with error estimator 
indicating how close they are approximating the high-dimensional data. 
RBM is a  relative recent approach to speed up the numerical simulation of parametric Partial Differential Equations (PDEs) 
\cite{NguyenVeroyPatera2005,Prudhomme_Rovas_Veroy_Maday_Patera_Turinici,Rozza_Huynh_Patera, CHMR_Sisc, ChenGottlieb}. 
It utilizes an Offline--Online computational decomposition strategy to produce surrogate solution (of dimension $N$) in a 
time that is of orders of magnitude shorter than what is needed by the underlying numerical solver of dimension ${\mathcal N} \gg N$ (called 
{\em truth} solver hereafter). 
The RBM relies on a projection onto a low dimensional space spanned by truth approximations at 
an optimally sampled set of parameter values \cite{
Almroth_Stern_Brogan,Fink_Rheinboldt_1,Noor_Peters,Porsching, Maday}. 
This low-dimensional manifold is generated by a greedy algorithm making use of a rigorous {\it a posteriori} error bounds
for the field variable and associated functional outputs of interest which also guarantees the fidelity 
of the surrogate solution in approximating the truth approximation. 

The RBD method acts in a similar fashion. Given the data matrix $X$ as in \eqref{eq:highd_x}, 
it iteratively builds up $Y$ \eqref{eq:lowd_y} whose column space approximates that of $X$. It starts with a randomly 
selected column of $X$ (or a user input if existent). At each step where we have $k$ vectors $\{y_1, \dots, y_k\}$, the next vector $y_{k+1}$ is found by 
scanning the columns of $X$ and locating the one whose error of projection into the current space ${\rm span}\{y_1, \dots, y_k\}$ 
is the largest. This process is continued until the maximum projection/compression error is small enough or until 
the limit on the size of the reduced space is reached. An important feature is an offline-online decomposition that allows the computation 
of the compression error, and thus the cost of locating $y_{k+1}$, to be independent of (the potentially large) $m$.

This paper is organized as follows. In Section \ref{sec:background}, we review the background material, mainly the RBM.  
Section \ref{sec:RBC} describes the reduced basis decomposition algorithm and discuss its properties. 
Numerical validations are presented in Section \ref{sec:numerical}, and finally some concluding remarks are offered in Section \ref{sec:conclusion}.

\section{Background}

\label{sec:background}

The reduced basis method 
was developed for use with finite element methods to numerically solve PDEs.
We assume, for simplicity, that the problems (usually parametric partial differential equations (PDE)) 
to simulate are written in the weak form:  find $u(\mu)$ in an Hilbert space $X$ such that
$a(u(\mu),v;\mu) = f(v;\mu), \,\,\,\, \forall v \in X$
where $\mu$ is an input parameter.  These simulations need to be performed  for many values of $\mu$  chosen
in a given parameter set ${\mathcal D}$. In this problem $a$ and $f$ are bilinear and linear forms,
respectively, associated to the PDE (with $a^\calN$ and $f^\calN$ denoting their numerical counterparts). 
 We assume that there is a numerical method to solve this problem and the solution $u^{\mathcal N}$, called 
 the ``truth approximation'' or ``snapshot'', is accurate enough for all $\mu \in {\mathcal D}$.

The fundamental observation utilized by RBM 
is that the parameter dependent solution 
$u^{\mathcal N}(\mu)$ is not simply 
an arbitrary member of the infinite-dimensional space associated with the PDE. 
Instead, the solution manifold ${\mathcal M} = \{u^{\mathcal N}(\mu),\,\mu \in {\mathcal
D}\}$ can typically be well approximated by a low-dimensional vector space. 
The idea is then to propose an approximation  of ${\mathcal M}$ by $W^N={\rm
span}\{u^{\mathcal N}(\mu_1),\,\dots,\,u^{\mathcal N}(\mu_N)\}$ where, $u^{\mathcal
N}(\mu_1),\,\dots,\,u^{\mathcal N}(\mu_N)$ are $N$ $(\ll\mathcal N)$ pre-computed truth approximations corresponding to
the parameters $\{\mu_1,\dots,\mu_N\}$ judiciously selected according to a sampling strategy \cite{Maday}. For
a given $\mu$, we now 
solve in $W^N$ for the reduced solution $u^{(N)}(\mu)$. The online computation is $\mathcal N$-independent, thanks to the
assumption that the (bi)linear forms are affine\footnote{$a(w,v;\mu) \,\equiv
\,\sum_{q=1}^{Q_a}\Theta_a^q(\mu)\,a^q(w,v),\quad \forall \,\,w,v \in X^\calN,$ $f(v;\mu) \,\equiv
\,\sum_{q=1}^{Q_f}\Theta_f^q(\mu)\,f^q(v),\quad \forall \,\,v \in X^\calN.$} and the fact that they can be
approximated by affine (bi)linear forms when they are nonaffine
\cite{Barrault_Nguyen_Maday_Patera,Grepl_Maday_Nguyen_Patera}. Hence, the online part is very efficient. In
order to be able to ``optimally'' find the $N$ parameters and to assure the fidelity of the reduced basis
solution $u^{(N)}(\mu)$ to approximate the truth solution $u^{\mathcal N}(\mu)$, we need an 
{\em a posteriori} error estimator $\Delta_N(\mu)$ 
which involves the residual $r(v,\mu) = f^\calN(v;\mu) - a^\calN (u^N(\mu),v;\mu)$ 
and stability information of the bilinear form
\cite{Machielis_Maday_Oliveira_Patera_Rovas,Maday_Patera_Rovas,Prudhomme_Rovas_Veroy_Maday_Patera_Turinici,
Rozza_Huynh_Patera,SenNatNorm}. 
With this estimator, we can describe briefly the classical {\bf greedy algorithm} used to find the $N$ parameters
$\mu_1, \dots, \mu_N$ and the space $W^N$.  
We first randomly
select one parameter value and compute the associated truth approximation. Next, we scan the entire 
(discrete) parameter space and for each parameter in this space compute its RB approximation $u^{(N=1)}$ and the error estimator
$\Delta_1(\mu)$. The next parameter value we select, $\mu_2$, is the one corresponding to the largest error
estimator. We then compute the truth approximation and thus have a new basis set
consisting 
of two elements.
This process is repeated until the maximum of the error estimators is sufficiently small.

The reduced basis method typically has exponential convergence with respect to the number of pre-computed
solutions
\cite{Maday_Patera_Turinici_2,BuffaMadayPateraPrudhommeTurinici2011,BinevCohenDahmenDevorePetrovaWojtaszczyk}.
This means that the number of pre-computed solutions can be small, thus the computational cost reduced
significantly, for the reduced basis solution to approximate the finite element solution reasonably well.
The author and his collaborators showed \cite{CHM_JCP}
that it works well even for a complicated geometric electromagnetic scattering problem 
that efficiently reveals a 
very sensitive angle dependence  (the object being stealthy with a particular configuration).

\section{Reduced basis decomposition}
\label{sec:RBC}

In this section, we detail our proposed methodology by stating the algorithm, studying the error evaluation, and pinpointing the 
computational cost.

\subsection{The algorithm}

\begin{algorithm}[h!]
  \caption{
    \hspace*{2mm}  
  Reduced Basis Decomposition \newline 
  \hspace*{2.9cm}  
  $(Y, T) = RBD(X, \epsilon_{\rm R}, d_{\rm max})$  
  }
  \label{alg:c_greedy}
  \begin{algorithmic}
\STATE {\bf 1.} Set $d = 1$, $E_{\rm cur} = +\infty$, and $i$ a random integer between $1$ and $n$.
\medskip
  \STATE {\bf 2.} 
  \WHILE { $d \le d_{\rm max}$ {\bf and} $E_{\rm cur} > \epsilon_{\rm R}$}
  \medskip
\STATE {\bf 2.1.} \framebox[0.85 \textwidth]{$v = X(:,i)$.}
\medskip
\STATE {\bf 2.2.} 
\framebox[0.85 \textwidth]
{\begin{minipage}{0.8 \textwidth}
Apply the modified Gram-Schmidt orthonormalization to obtain the $d^{\rm th}$ basis of the compressed space
\FOR{$j = 1:d-1$}
\STATE $v = v - (v \cdot \xi_j) \, \xi_j$.
\ENDFOR
\IF{$\lVert v \rVert < \epsilon_{\rm R}$}
\STATE $Y = Y(:,1:d-1)$
\STATE $T = T(:,1:d-1)$
\STATE Break;
\ELSE
\STATE $\xi_d = \frac{v}{\lVert v \rVert}$, $Y(:,d) = \xi_d$.
\STATE $T(d,:) = \xi_d' X$.
\ENDIF
\medskip
\end{minipage}
}
\medskip
\STATE {\bf 2.3.} 
\medskip
\framebox[0.85 \textwidth]{
\begin{minipage}{0.8 \textwidth}
$E_{\rm cur} = \displaystyle{\max_{j \in \{1, \dots, n\}} \lVert X(:,j) - Y(:,1:d) T(:,j)\rVert}$\\
and \\
$i = \displaystyle{\argmax_{j \in \{1, \dots, n\}} \lVert X(:,j) - Y(:,1:d) T(:,j)\rVert}$
\medskip
\end{minipage}
}
\STATE {\bf 2.4.} 
\framebox[0.85 \textwidth]{
\begin{minipage}{0.8 \textwidth}
\IF{$E_{\rm cur} \le \epsilon_R$}
\STATE $Y = Y(:,1:d)$, $T = T(1:d,:)$.
\ELSE
\STATE $d = d + 1$.
\ENDIF
\medskip
\end{minipage}
}
\medskip
\ENDWHILE
  \end{algorithmic}
\end{algorithm}
At the heart of the method stated in Algorithm \ref{alg:c_greedy} is a greedy algorithm similar to that used by RBM. 
It builds the reduced space dimension-by-dimension. At each step, the {\em greedy} decision for the best next dimension 
to pursue in the space corresponding to the data is made by examining an error indicator quantifying the discrepancy between the uncompressed data 
and the one compressed into the current (reduced) space. 

In the context of the RBM, we view each column (or row if we are compressing the row space) 
of the matrix as the fine solution of certain (virtual) parametric PDE 
with the (imaginary) parameter taking a particular value. Since this solution is {\em explicitly} given already by the data, the fact that 
the PDE and the parameter are absent does not matter. 
Once this {\em common mechanism} satisfied by each column (or row) is identified, 
the greedy algorithm still relies on an accurate and efficient estimate quantifying the error between the original data and the 
compressed one. This will be the topic of the next subsection. 

To state the algorithm, we assume that we are given a data matrix $X \in {\mathbb R}^{m \times n}$, the largest dimension $d_{\rm max} < n$ that the practitioner wants to retain, and 
a tolerance $\epsilon_{\rm R}$ capping the discrepancy between the original and the compressed data. 
The output is the set of bases for the compressed data (a low-dimensional approximation of the original data)
 $Y \in {\mathbb R}^{m \times d}$ and the transformation matrix $T \in {\mathbb R}^{d \times n}$. Here, 
 $d \le d_{\rm max}$ is the actual dimension of the compressed data. 
 
 With this output, we can
\begin{description} 
\item [{\bf Compress.}] We represent any data entry $X(:,j)$, the $j^{\rm th}$ column of $X \in {\mathbb R}^m$, by the $j^{\rm th}$ column of $T$, 
$T(:,j) \in {\mathbb R}^d$, with usually $d \ll m$.
\item [{\bf Uncompress.}] An {\em approximation} of the data is reconstructed by
$$X(:,j) = Y\,T(:,j).$$
\item [{\bf Evaluate the compression of out-of-sample data.}] Given any $v \in {\mathbb R}^{m \times 1}$ that is not 
 equal to any column of $X$, its compressed representation in ${\mathbb R}^{d \times 1}$ is 
 $$v_{\rm C} = Y' v.$$
\end{description}

\subsection{Efficient quantification of the error}

A critical part to facilitate the greedy algorithm and make the algorithm realistic is an efficient mechanism 
measuring (or estimating) the error $v - v_{\rm C}$ under certain norm, $\lVert v - v_{\rm C}\rVert$, in Step $2.3$ of the algorithm. 
In this work, we are using the $A-$norm defined as follows. 
For a given symmetric and positive definite matrix $A \in {\mathbb R}^{m \times m}$, 
the $A-$norm of a vector $v \in {\mathbb R}^{m \times 1}$ is defined by
$$\lVert v \rVert_A := \sqrt{v' A v}.$$

For $v$ being any column of the data matrix $X$ and $v_{\rm C}$ its low-dimensional approximation $v_{\rm C} = Y \vec{c}$, it is easy to see that
\begin{alignat}{1}
\label{eq:evalerr}
\lVert v - v_{\rm C}\rVert_A^2 & = v' A v - 2 v_{\rm C}' A v + v_{\rm C}' A v_{\rm C}\\
\nonumber
& = v' A v - 2 \vec{c}\,' Y' A v + \vec{c}\,' Y' A Y \vec{c}.
\end{alignat}

The choice of $A$ reflects the criteria of the data compression. Typical examples are:
\begin{itemize}
\item [1.] {\bf Identity: } Equal weights are assigned to each component of the data entry. 
This makes the quality of compression uniform. 
In this case, the evaluation of \eqref{eq:evalerr} is greatly simplified and the algorithm is 
the fastest as shown below by the numerical results. 
\item [2.] {\bf General diagonal matrix: } This setting can be used if part of each data entry 
needs to be preserved better and other parts can afford less fidelity. 
\item [3.] {\bf General SPD matrix: } This most general case can be helpful if the goal is to preserve 
data across different entries anisotropiclly.
\end{itemize}
The goal is then to evaluate the error through \eqref{eq:evalerr} as efficiently as possible for any given $\vec{c}$. This is achieved by 
employing an offline-online decomposition strategy where the $\vec{c}$-independent parts are evaluated beforehand (offline) enabling a quick 
turnaround time for any given $\vec{c}$ encountered online. The specifics are given in the next subsection.

\subsection{Computational cost and implementation aspects}

The Offline-Online decomposition of the computations and their complexities are as follows. Here, we use 
$nnz(A)$ to denote the number of nonzero entries of a sparse matrix $A$.
\begin{description}
\item[{\bf Offline}] The total cost is of order $$\left(m + nnz(A)\right) \, d_{\rm max}^2 + \left(d^3_{\rm max} + nnz(A)\right) n.$$
\begin{description}
\item[{\bf Offline MGS}] Every basis needs to orthogonalized against the current set of bases. 
The total cost is of order $m \, d_{\rm max}^2$.
\item[{\bf Offline Calculation of Errors}] The next basis is located by comparing each column with its compressed version into 
the current space. To enable that, we encounter the following computational cost:
\begin{description}
\item [{\bf Pre-computation}] of $diag(X'AX)$ (for $v' A v$ in \eqref{eq:evalerr}) and $AX$ (for $Av$ in \eqref{eq:evalerr}). The cost is of order $nnz(A) \, n$.
\item [{\bf Expansion}] of $Y'AX$ and $Y' A Y$. The former takes time of order $nnz(A) \, d_{\rm max} \, m$, 
and the latter of order $nnz(A) \, d^2_{\rm max}$.
\end{description}
\item[{\bf Offline Searching}] After these calculations, 
the comparison between the original and compressed data is then only dependent on the size of $\vec{c}$ (which 
is also the number of columns for $Y$). The complexity is of order $d^2_{\rm max}$. It will be repeated for up to $n$ times in the searching process 
of step 2.3 of the algorithm for each of the up to $d_{\rm max}$ basis elements. The total cost is at the level of $n\,d^3_{\rm max}$.
\end{description}
\item[{\bf Online}] Given any (possibly out-of-sample) data $v \in {\mathbb R}^{m \times 1}$, its coefficients in the 
compressed space is obtained by evaluating $\vec{c} = Y' v$. The cost is of order $m \, d_{\rm max}$. The decoding 
($Y\,\vec{c}$) can be done with the same cost. 
The online computation has complexity of order $$m \, d_{\rm max}.$$
We remark that, if the actual practice does not requires forming $v_{\rm C}$ (e.g. clustering and classification etc) and so we only work with 
the coordinates $\vec{c}$ of $v$ in the compressed space, then the online cost will be independent of $m$ and thus much smaller. 
\end{description}

\section{Numerical Results}

In this section, we test the reduced basis decomposition on image compression, and data compression.
Lastly, we devise a simple face recognition algorithm based on RBD and test it on a database of $575$ images 
while comparing RBD with $6$ other face recognition algorithms. 
The computation is done, and thus the speedup numbers reported herein should be understood as, in 
Matlab 2014a on a 2011 IMac with a $3.4$ GHz Intel Core i7 processor.

\label{sec:numerical}

\subsection{Image Compression and comparison with SVD}
\begin{figure}[htbp]
\begin{center}
\includegraphics[width=0.4\textwidth]{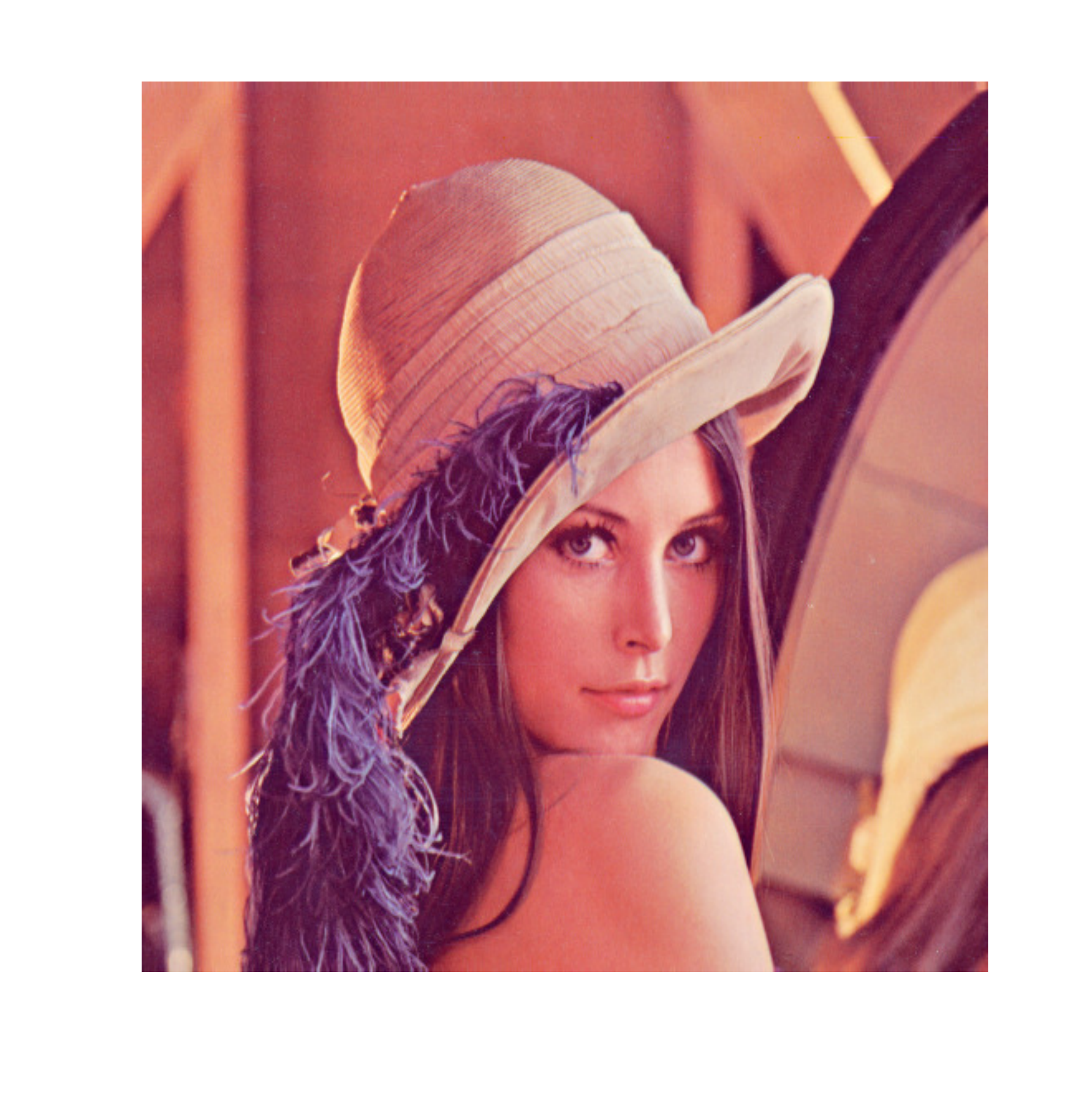}
\includegraphics[width=0.4\textwidth]{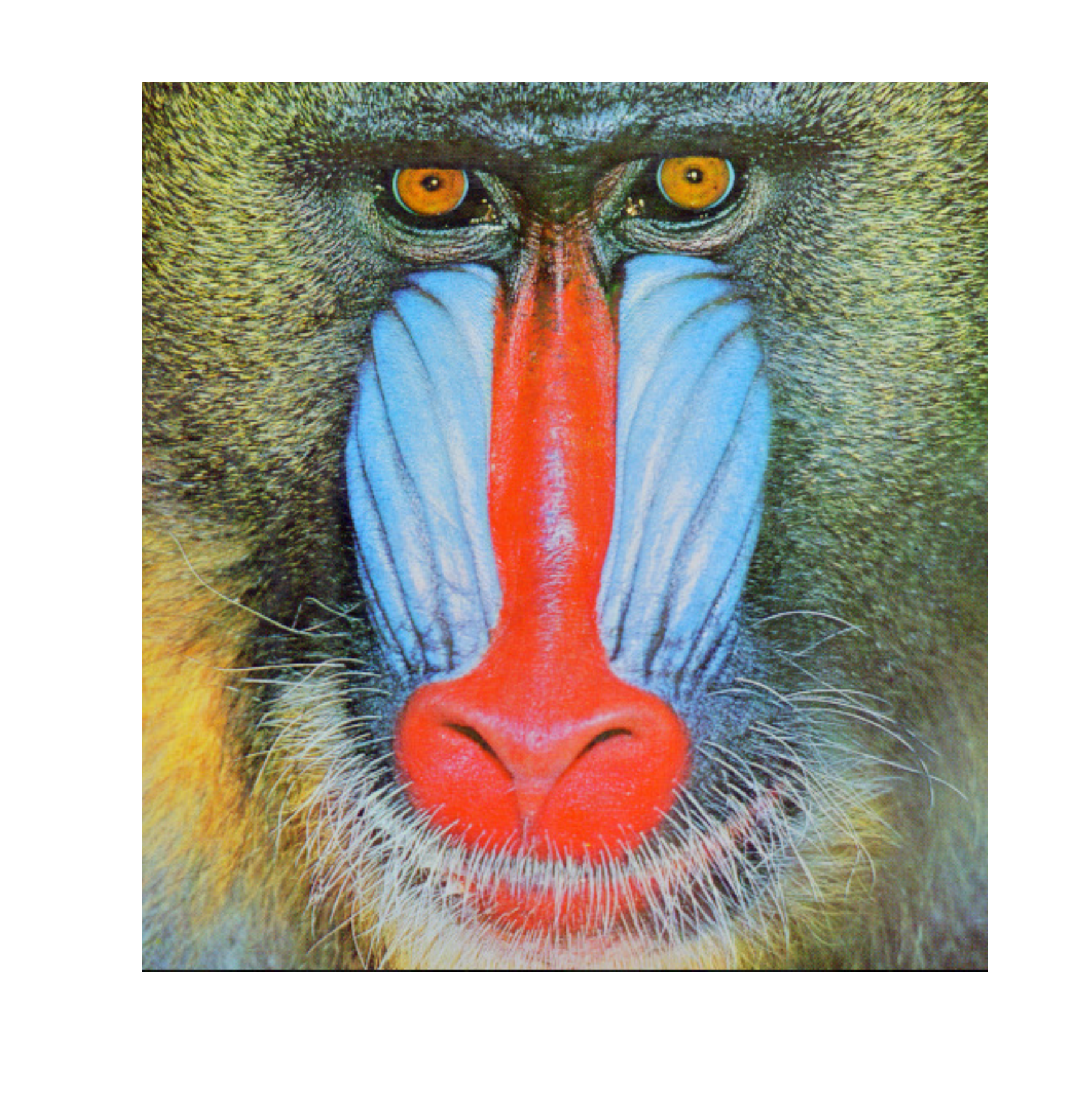}
\caption{Original pictures: Lena and Mandrill.}
\label{fig:original}
\end{center}
\end{figure}

\begin{figure}[htbp]
\begin{center}
\includegraphics[width=0.32\textwidth]{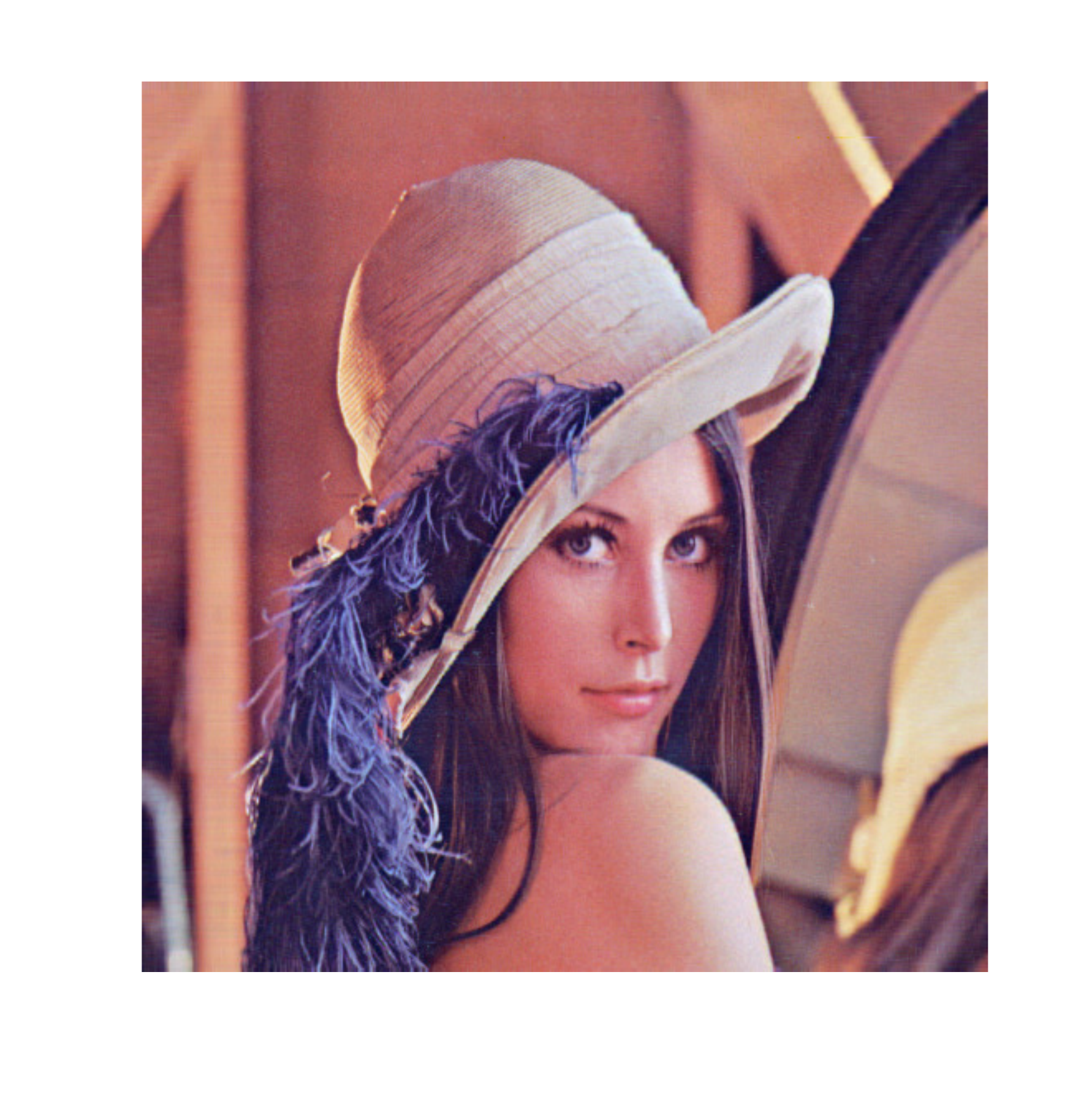}
\includegraphics[width=0.32\textwidth]{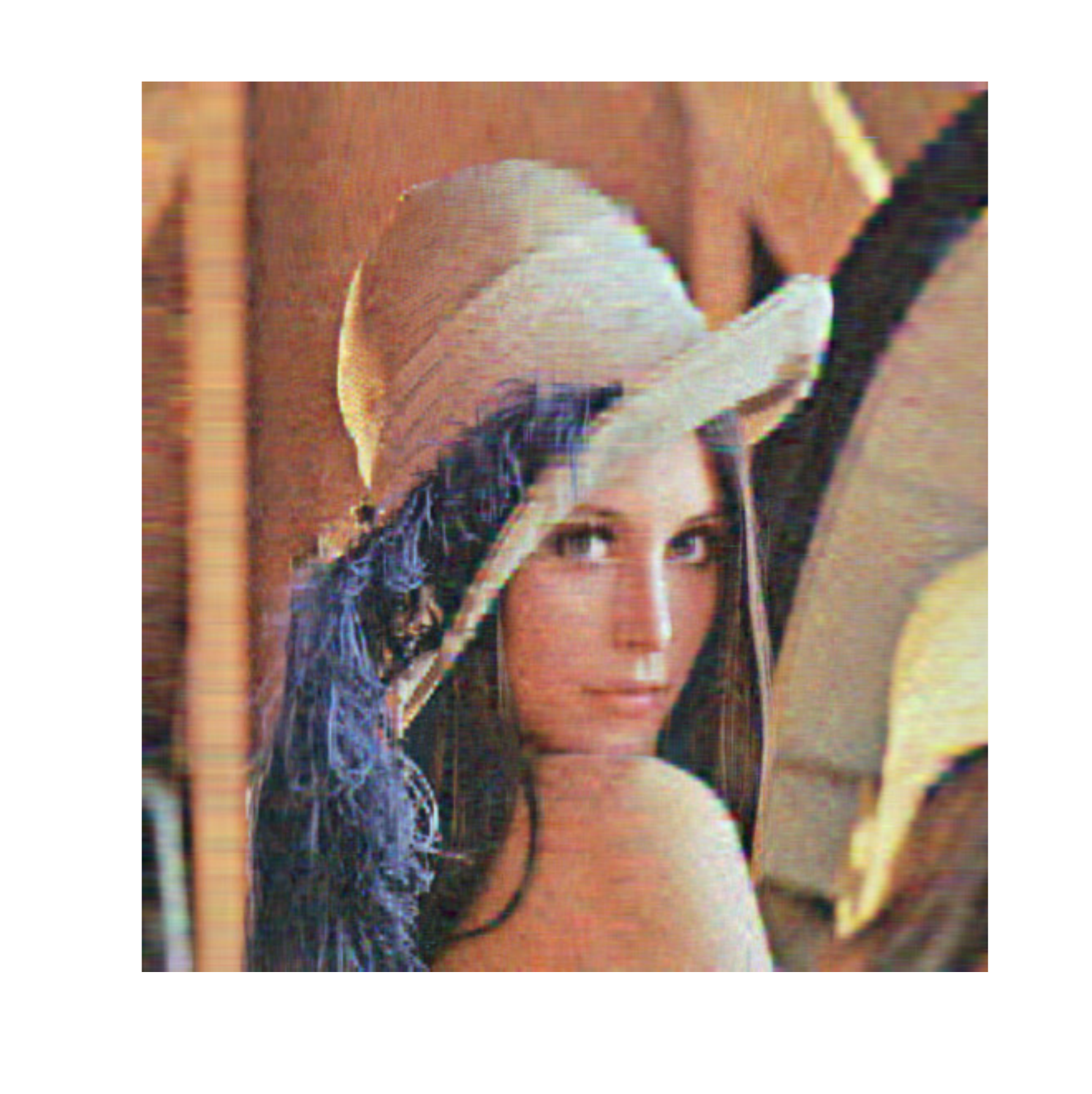}
\includegraphics[width=0.32\textwidth]{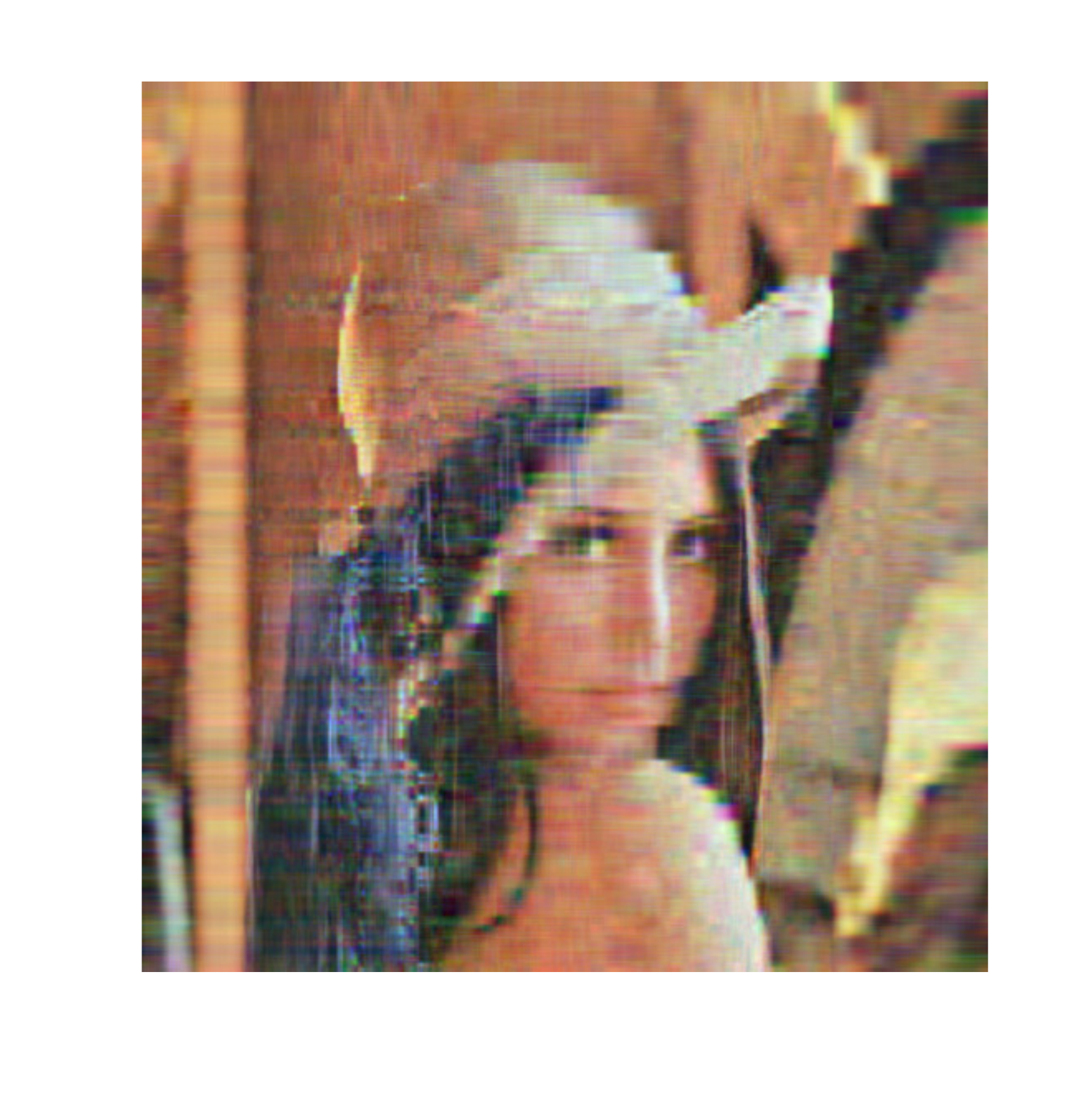}\\
\includegraphics[width=0.32\textwidth]{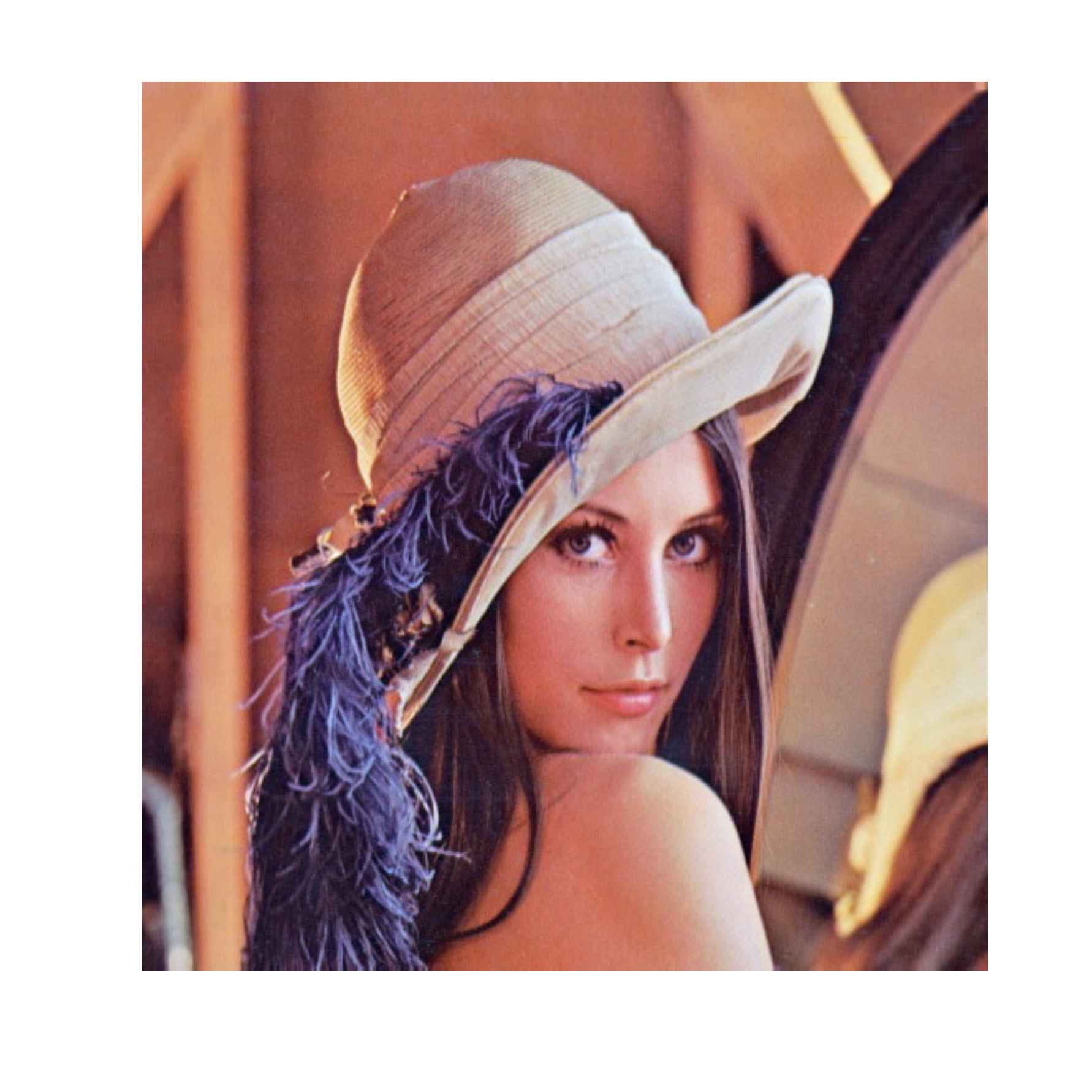}
\includegraphics[width=0.32\textwidth]{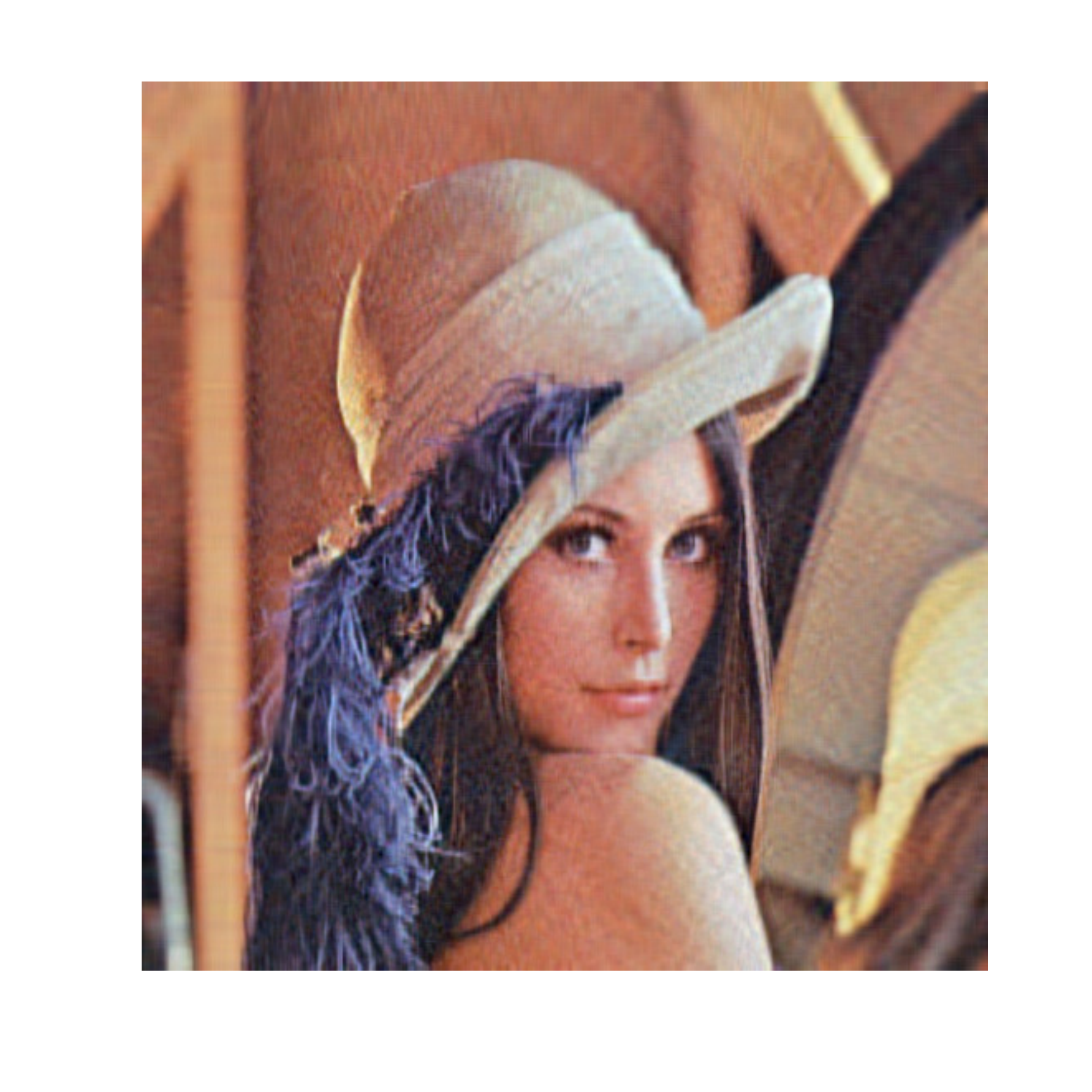}
\includegraphics[width=0.32\textwidth]{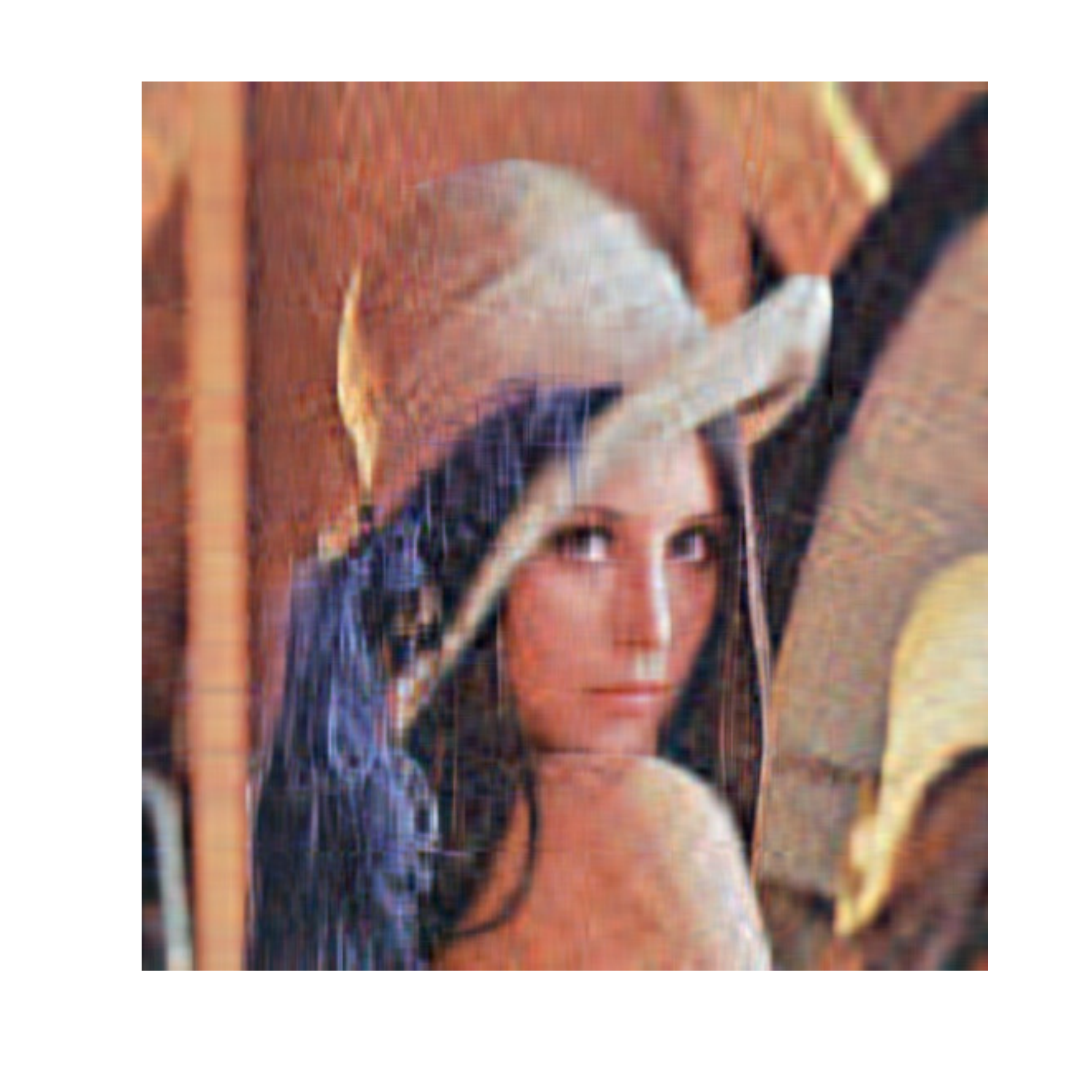}\\
\includegraphics[width=0.32\textwidth]{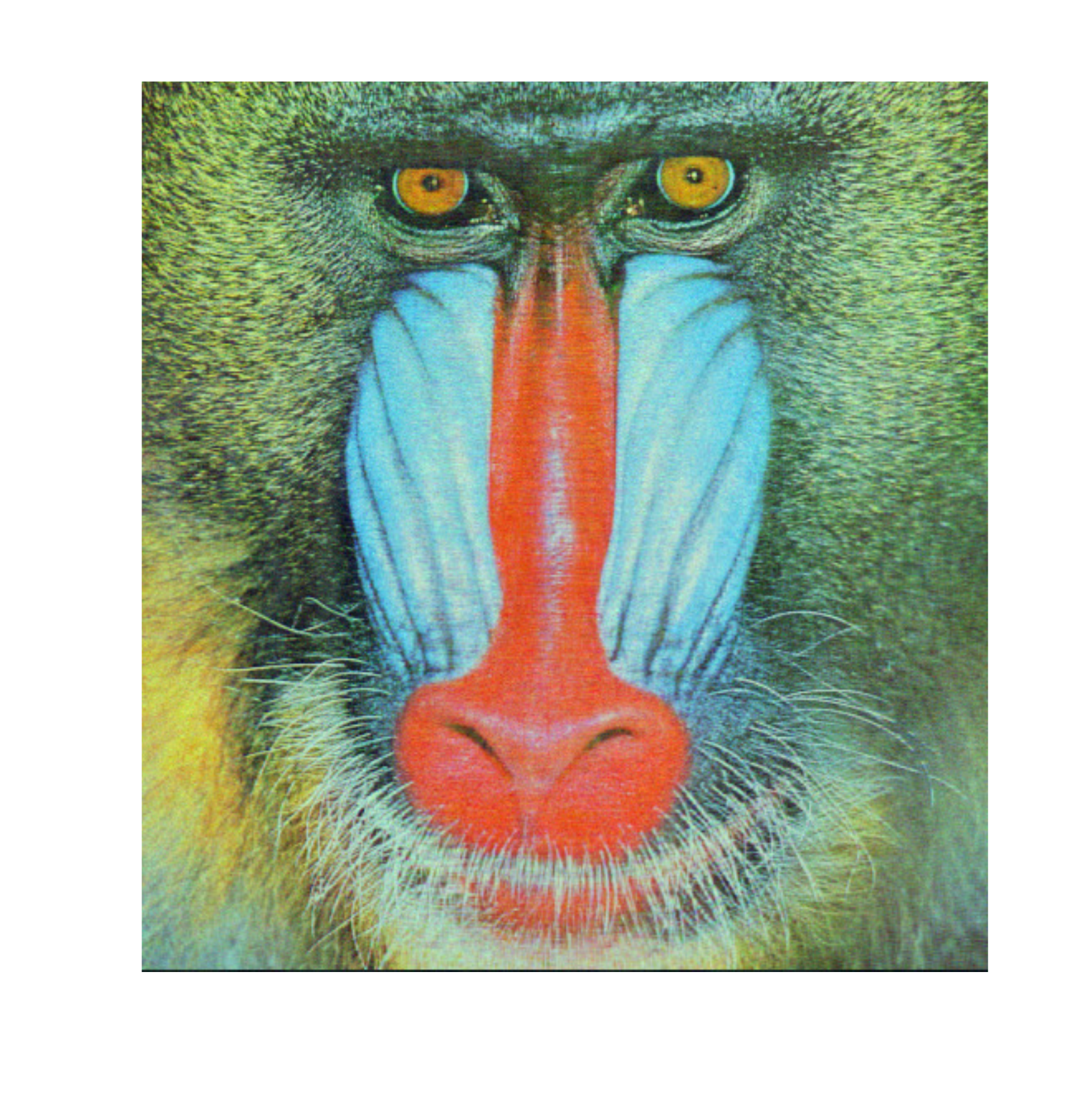}
\includegraphics[width=0.32\textwidth]{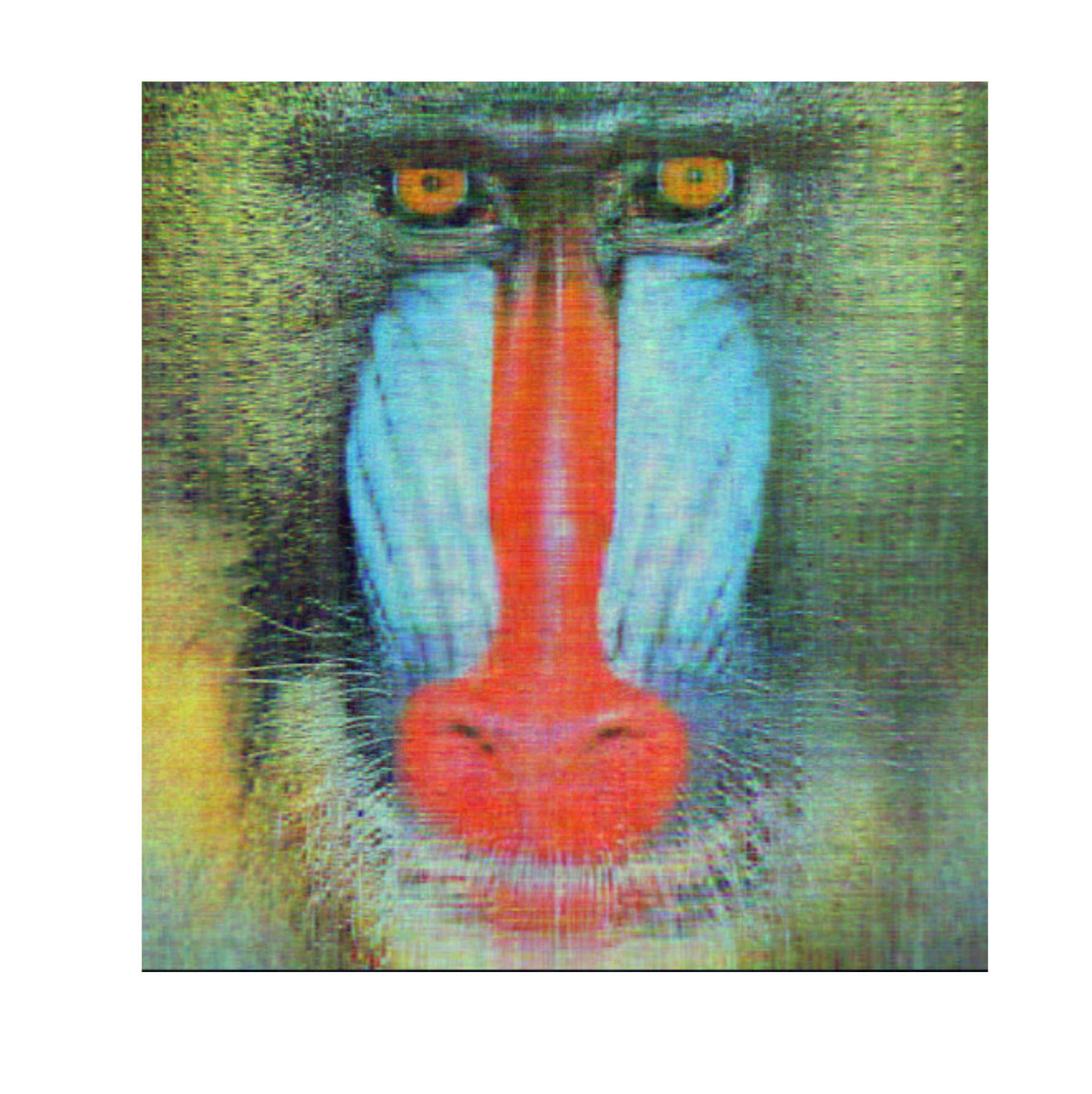}
\includegraphics[width=0.32\textwidth]{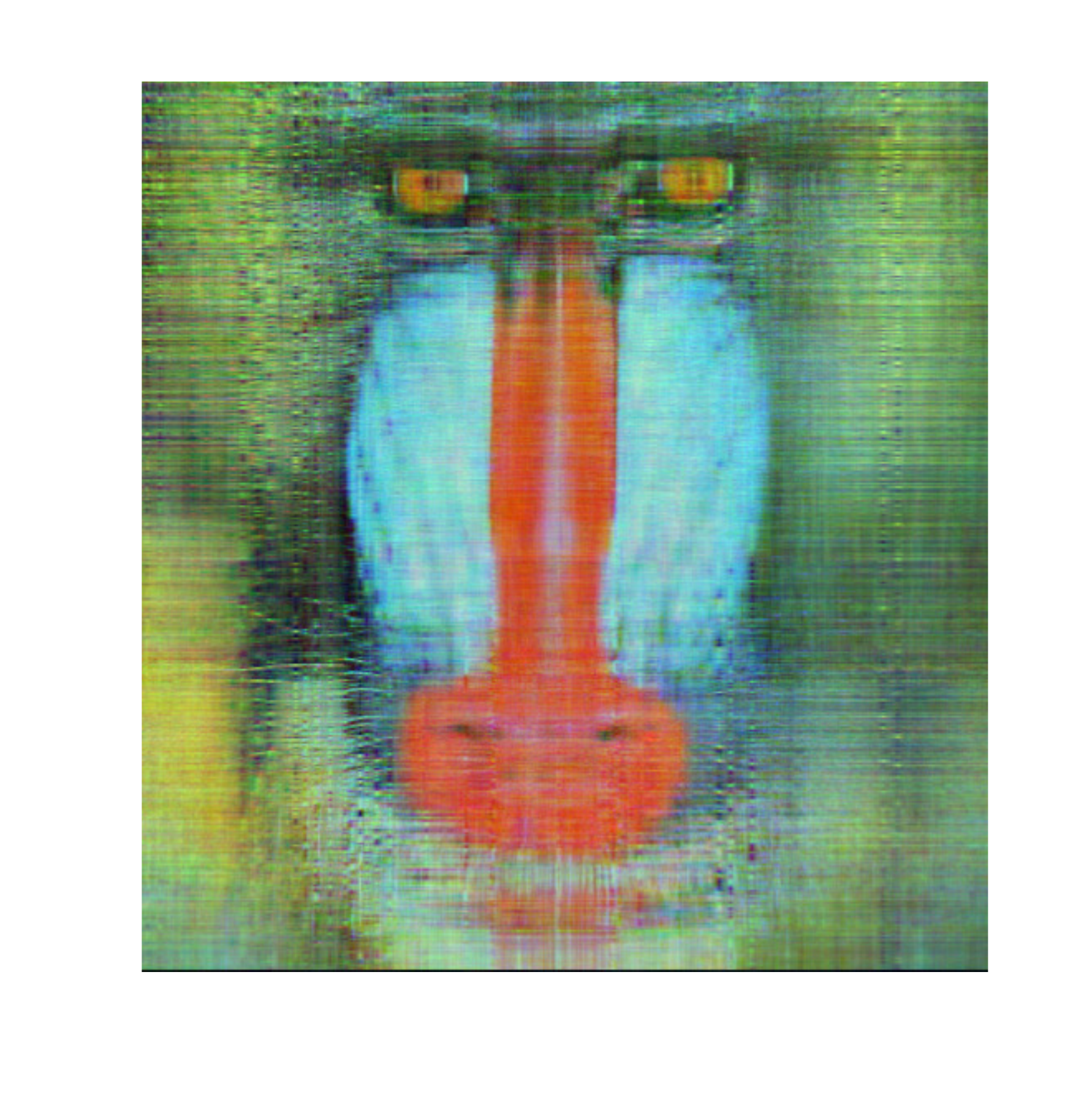}\\
\includegraphics[width=0.32\textwidth]{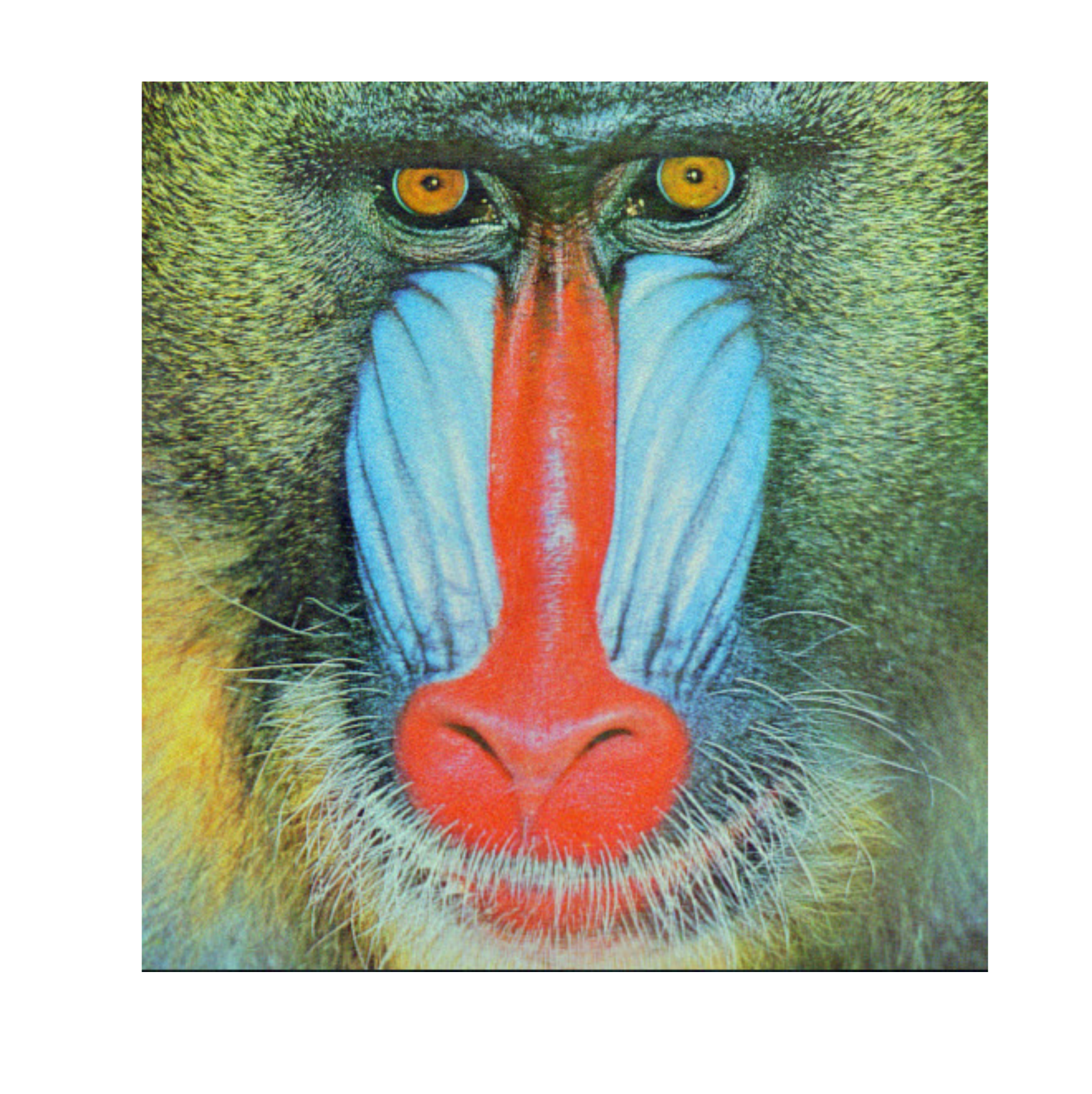}
\includegraphics[width=0.32\textwidth]{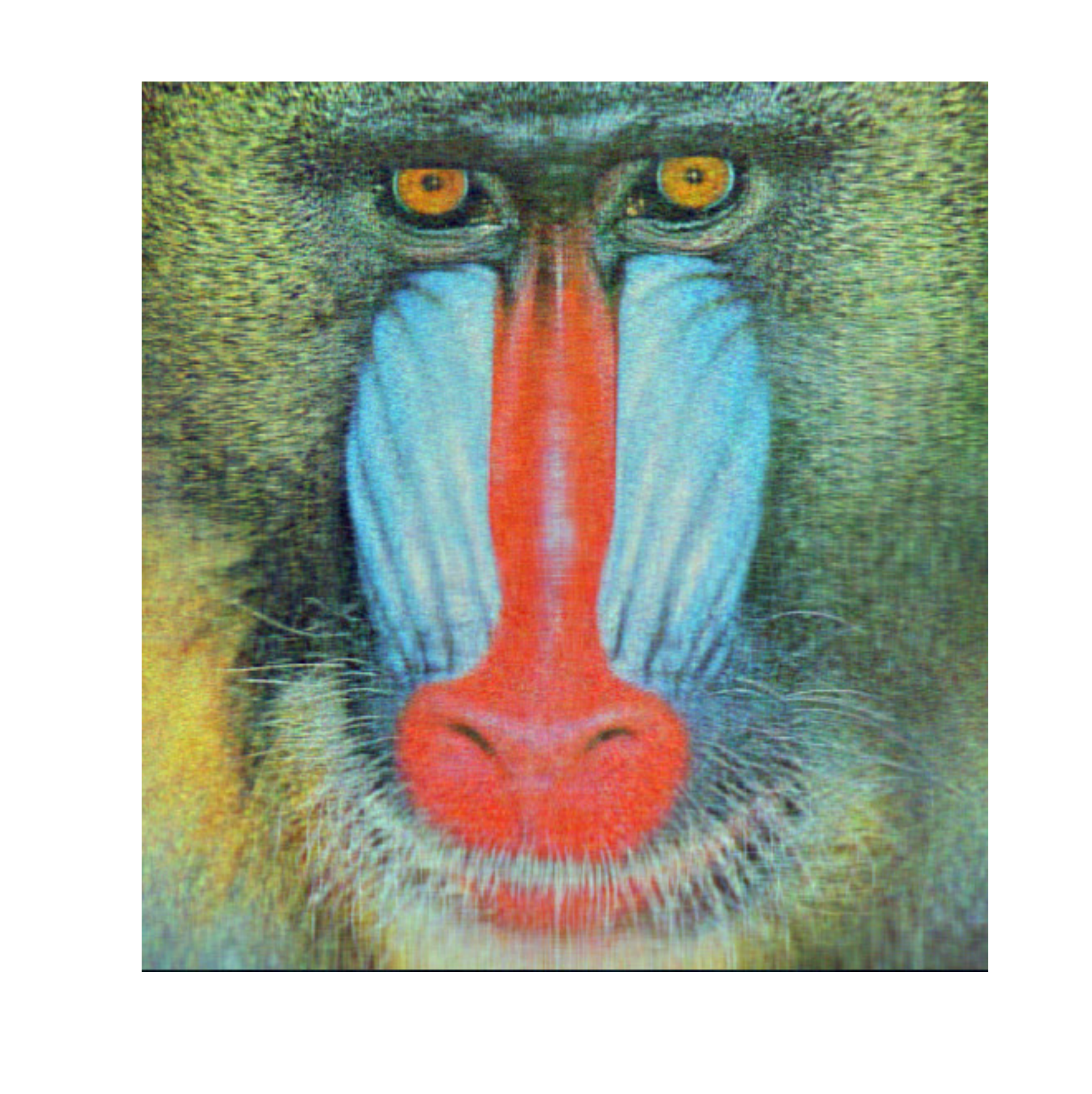}
\includegraphics[width=0.32\textwidth]{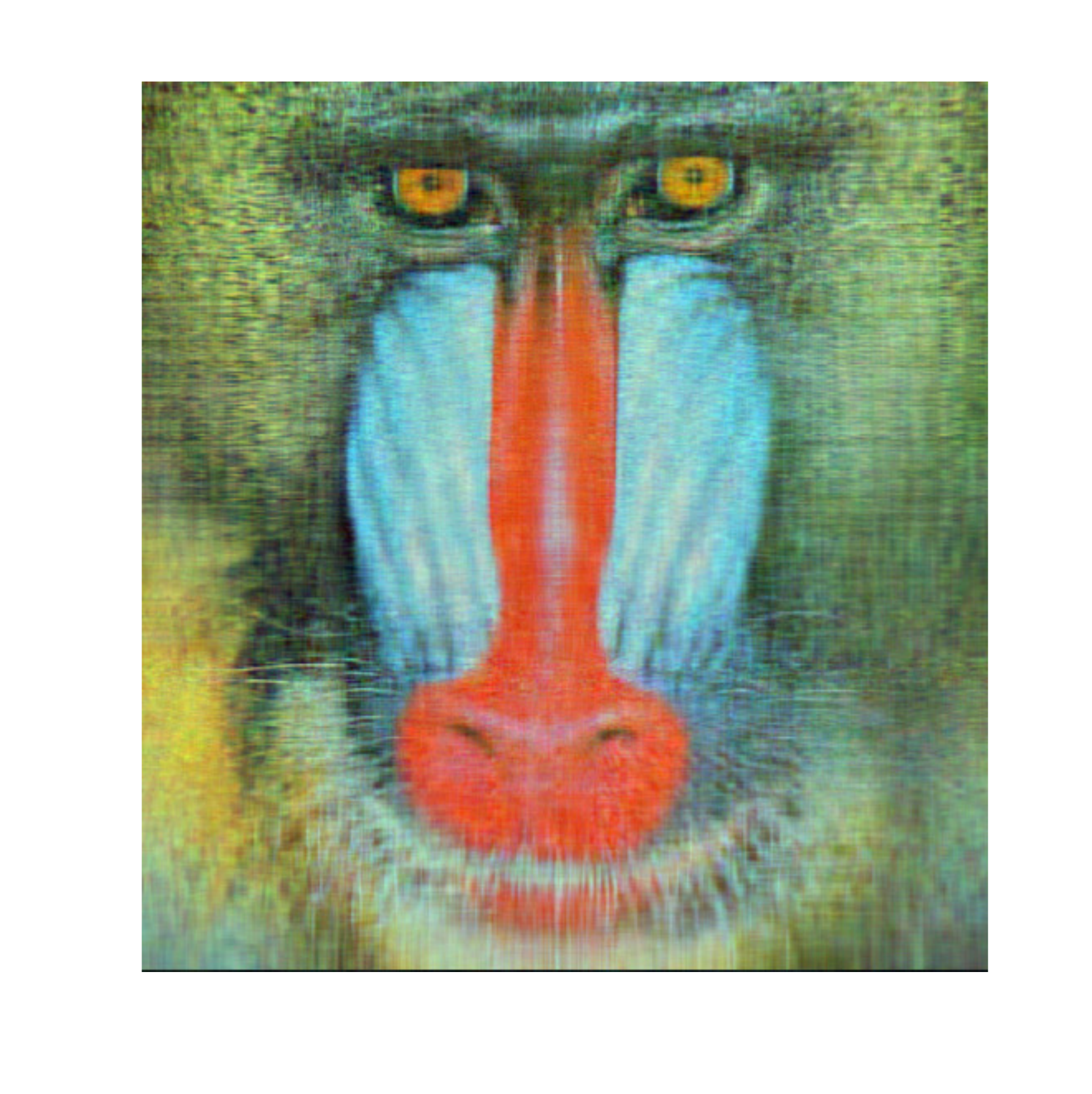}
\caption{Lena and Mandrill compressed. $d = 170, 51, 25$ from left to right. The first and third row are 
from Reduced Basis Decomposition, and the second and fourth are from Singular Value Decomposition.}
\label{fig:LenaResult}
\end{center}
\end{figure}

We first test it on compressing two standard images Lena and Mandrill in Figure \ref{fig:original}. They both have 
an original resolution of $512 \times 512$. 
We take $A = I$ and test the algorithm. For each component of every image, we run the algorithm with 
$d_{\rm max} \in \{170, 51, 25\}$ which implies a compression ratio of $33\%$, $10\%$, and $5\%$ respectively. 
The resulting images (formed by multiplying the corresponding $Y$ and $T$ together) are shown  on the $1^{st}$ and $3^{rd}$ row of 
Figures \ref{fig:LenaResult}. As a comparison, we run SVD and obtain the reconstructed 
matrices with the first $d_{\rm max}$ singular values accordingly. The resulting images are on the second and last row. Clearly, SVD provides the best quality 
pictures among all possible algorithms (and thus better than what RBD provides). 
However, we see that the RBD pictures are only slightly blurrier. Moreover, it takes much less 
time. In fact, we show the comparison in time between SVD and RBD in Table \ref{tab:LenaTime}. 
We see that, when $d = 51$, RBD is three times faster than SVD and seven times faster when $d = 25$. Here the SVD time 
is the shorter between those taken by {\verb svd } and {\verb svds } commands in Matlab.

\begin{table}[htbp]
\begin{center}
\renewcommand{\arraystretch}{1.3}
\begin{tabular}{|c|c|c|}
\hline
Picture & \quad RBD \quad & \quad SVD \quad \\
\hline
Lena, $d = 170$ & 3.57 & 1 \\
Lena, $d = 51$ & 0.30 & 1 \\
Lena, $d = 25$ & 0.14 & 1 \\
\hline
\quad Mandrill $d = 170$ \quad & 3.22 & 1 \\
Mandrill $d = 51$ & 0.31 & 1 \\
Mandrill $d = 25$ & 0.14 & 1 \\
\hline
\end{tabular}
\caption{Relative computational time for image compression.}
\label{tab:LenaTime}
\end{center}
\end{table}

\subsection{Data Compression}
Here, we test the algorithm on a few artificially-generated data sets. Given a function $f(x, y)$, the data 
denoted by $f(\calD)$ is constructed by evaluating $f$ on a 
uniform tensorial grid $\calD := (x_i, y_j)_{i,j = 1}^n$.

\subsubsection{Exact reconstruction}
For tensorial functions such as those listed in Table \ref{tab:rank123funcs} 
with their corresponding $d$ values, the RBD method detects the optimal dimension, 
stops the greedy algorithm after $d$ steps and decompose the matrix $f(\calD)$ accordingly, 
that is, as an {\em exact} product of $n \times d$ and $d \times n$ matrices.
\begin{table}
\begin{center}
\renewcommand{\arraystretch}{1.5}
\begin{tabular}{|c|c|}
\hline
Function $f(x,y)$ & Intrinsic dimension \\
\hline
$\sin(\pi x) \cos(\pi y)$ & 1 \\
\hline
$\sin(\pi x) \cos(\pi y) + 0.1 \sin(10 \pi x) \cos(10 \pi y)$ & 2 \\
\hline
\begin{minipage}{3 in}
\begin{center}
$\sin(\pi x) \cos(\pi y) + 0.1 \sin(10 \pi x) \cos(10 \pi y) $ 
$+ 0.01 \sin(100 \pi x) \cos(100 \pi y)$ 
\end{center}
\end{minipage}
& 3 \\
\hline
\end{tabular}
\caption{Three functions with low intrinsic dimensions that can be compressed by RBD exactly.}
\label{tab:rank123funcs}
\end{center}
\end{table}

\subsubsection{Approximate reconstruction}

\begin{figure}[htbp]
\begin{center}
\includegraphics[width=0.49\textwidth]{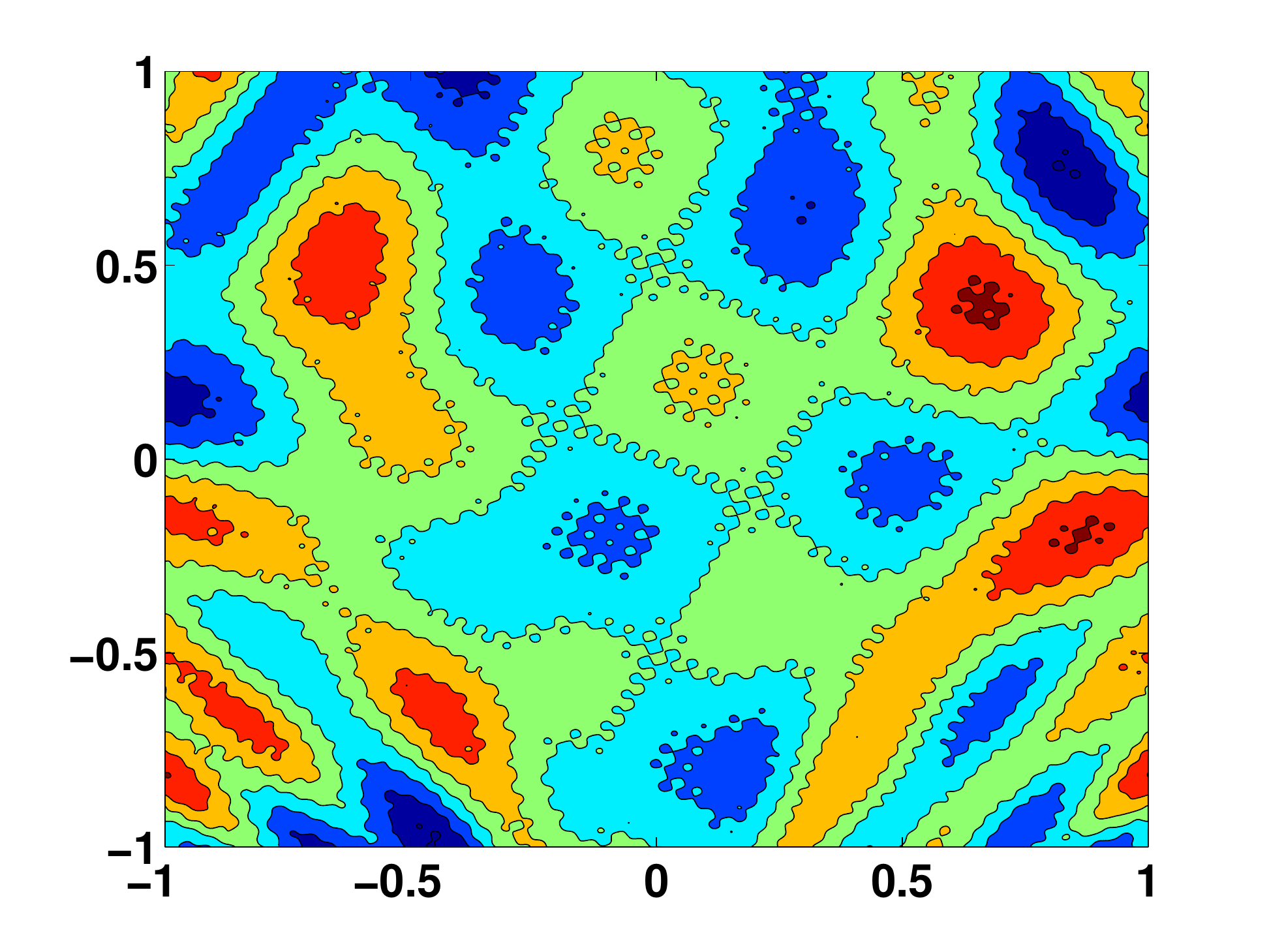}
\includegraphics[width=0.49\textwidth]{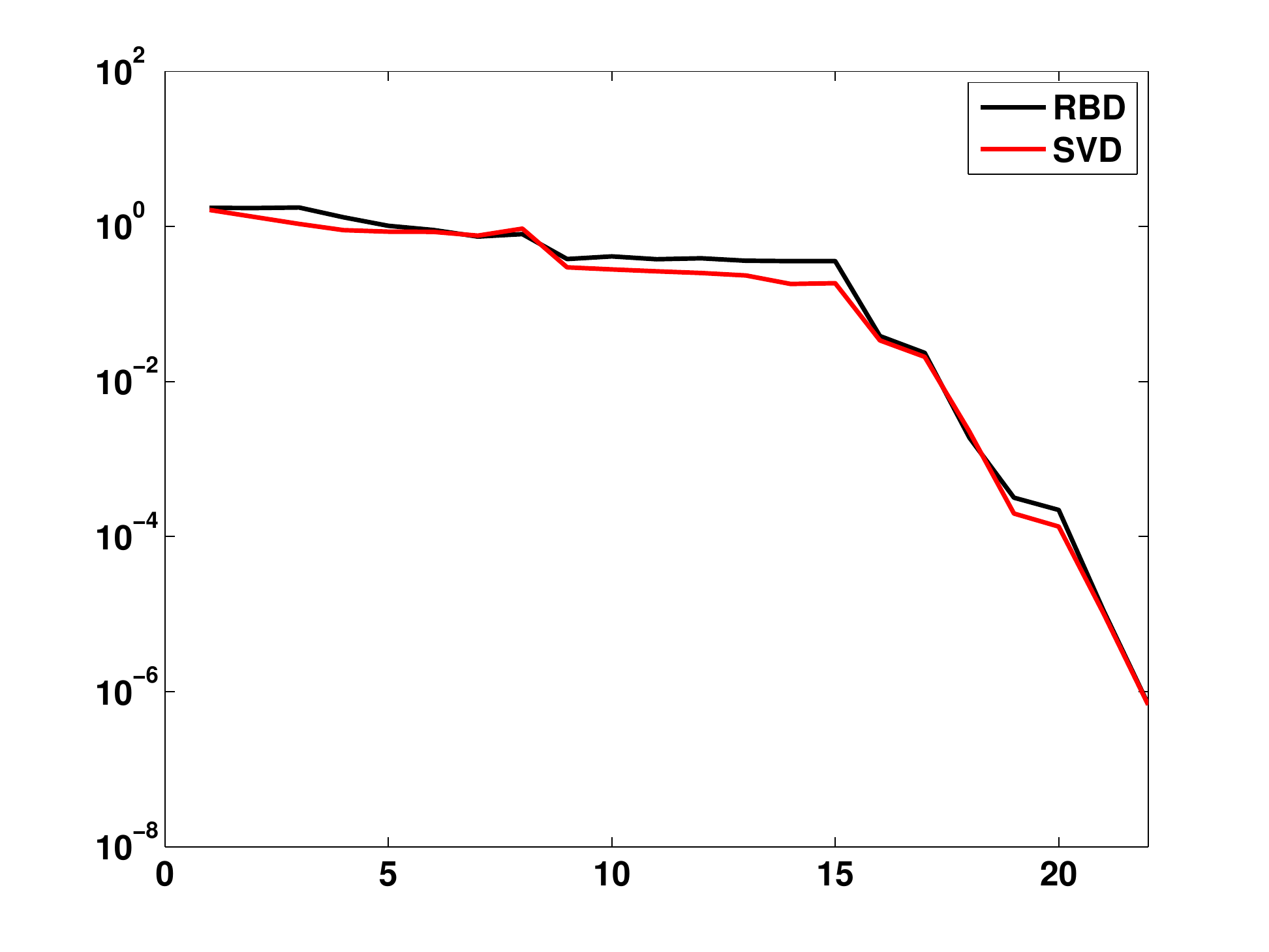}
\caption{Artificial dataset: The reconstructed contour plot from compressed data (left), and the comparison of the history of convergence (RBD vs SVD) as $d$ increases (right).}
\label{fig:FakeDataResult}
\end{center}
\end{figure}

Here, we 
set $f(x, y) = 0.6 f_1(x, y) + 0.1 f_2(x, y) + 0.01 f_3(x, y)$ with:
\begin{align*}
f_1(x, y) & = \sin(\pi (x + 2 y)) \cos(\pi (2 x - y)) \\
f_2(x, y) & = \sin(10 \pi (x - 3 y)) \cos(10 \pi (3 x + y))\\
f_3(x, y) & = \sin(3 \pi x^2 y) \cos(6 \pi \frac{\sqrt{|x|}}{y + 2}).
\end{align*}
and let $\calD$ be a $5001 \times 5001$ uniform grid on $[-1, 1] \times [-1,1]$. 
Setting $d = 22$, the method extracts $22$ columns and decompose $f(\calD)$ by a product of two matrices of size $5001 \times 22$ and $22 \times 5001$.  
The compression ratio is larger than $110$. 
More importantly, the reconstruction plotted in Figure \ref{fig:FakeDataResult} Left, has point-wise error below $10^{-6}$. 

We calculate the point-wise reconstruction error for reduced basis decomposition $e_R(d) = \lVert f(\calD) - Y(:,1:d) T(1:d,:)\rVert$. As a comparison, we 
calculate the first $22$ singular values $s_i$ of $f(\calD)$, the corresponding singular vectors $(u_i, v_i)$, and the reconstruction error 
$e_S(d) = \lVert f(\calD) - \sum_{i = 1}^d s_i u_i v_i'\rVert$. 
These two errors are 
plotted in  Figure \ref{fig:FakeDataResult} Right. We see that RBD {\em matches} SVD in terms of accuracy. We emphasize that what is striking 
is its efficiency. 
The RBD code, as implemented by the author\footnote{www.faculty.umassd.edu/yanlai.chen/RBD} is $16$ times faster than the 
{\verb svds }
command in Matlab.

\subsection{Face Recognition}

Here, we demonstrate the superior efficiency and accuracy of the RBD method on a classical classification task -- face recognition. 
The goal of face recognition is to recognize subjects based on facial images. 
It has important applications in areas ranging from surveillance, authentication, to human-computer interaction etc.

\begin{table}
\begin{center}
\renewcommand{\arraystretch}{1.3}
\begin{tabular}{|c|c|c|}
\hline
Data set & No of classes & No of samples per class \\
\hline
UMIST & 20 & 19--48 \\
\hline
\end{tabular}
\caption{Data set information.}
\label{tab:datasetinfo}
\end{center}
\end{table}
\begin{figure}[htbp]
\begin{center}
\includegraphics[width=\textwidth]{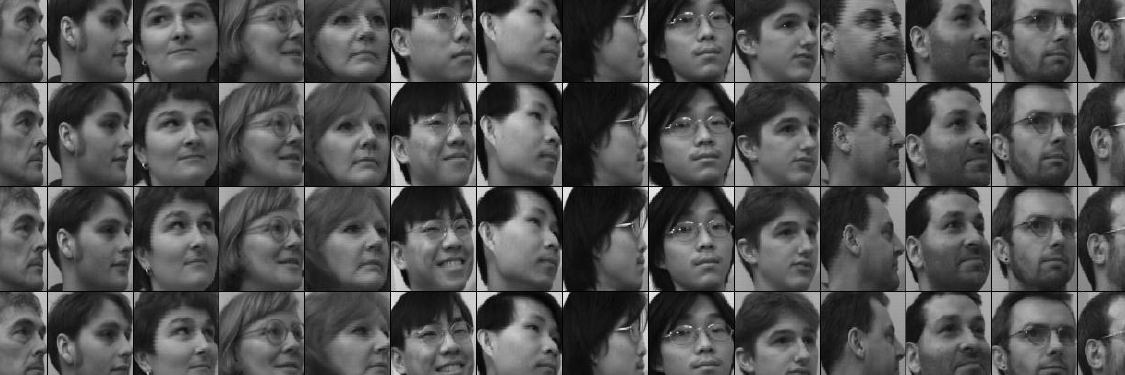}
\caption{A snapshot of the UMIST data set.}
\label{fig:UMIST}
\end{center}
\end{figure}
\begin{figure}[htbp]
\begin{center}
\includegraphics[width=0.32\textwidth]{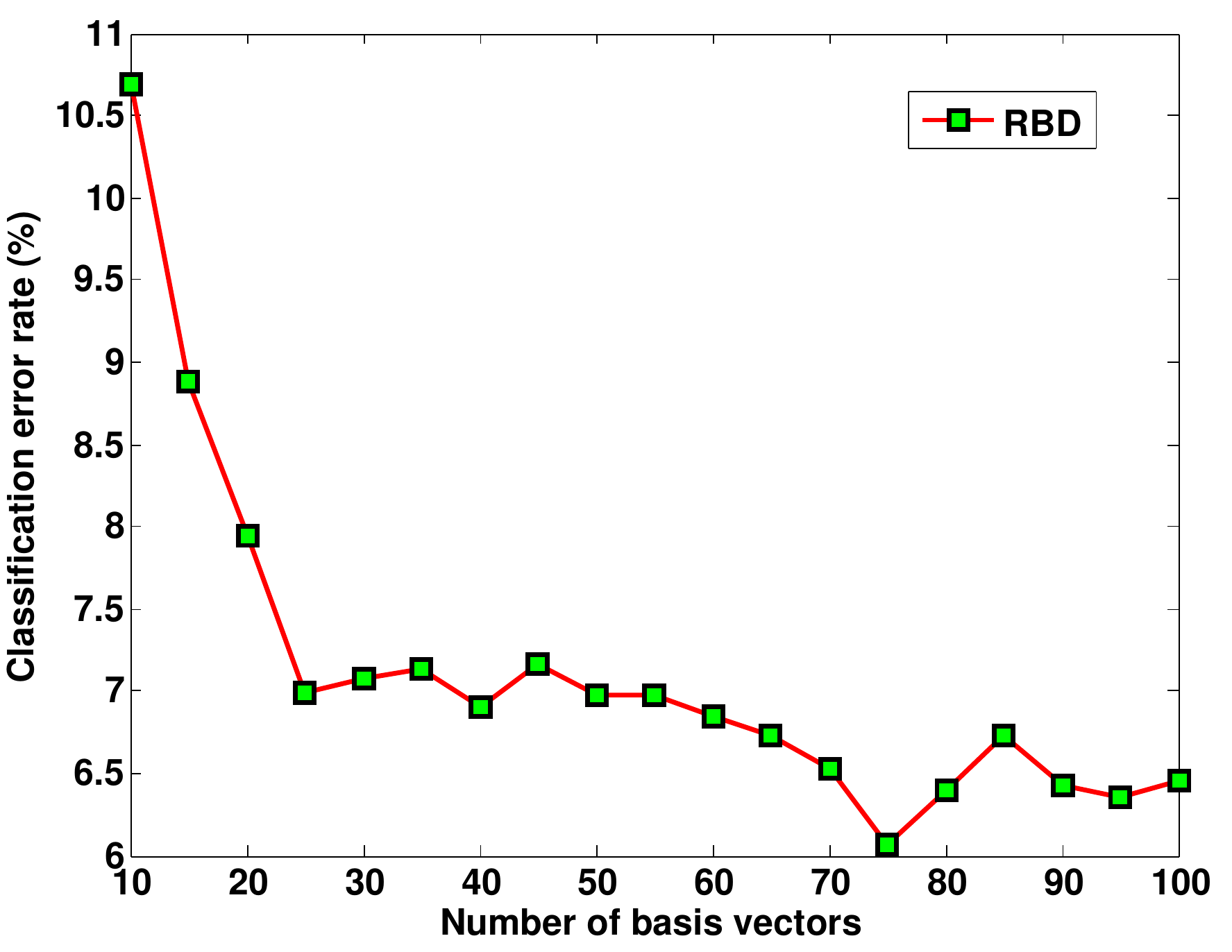}
\includegraphics[width=0.32\textwidth]{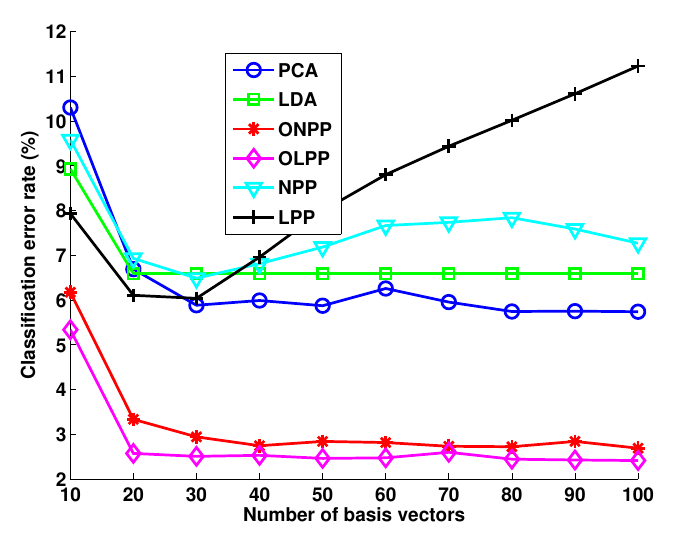}
\includegraphics[width=0.32\textwidth]{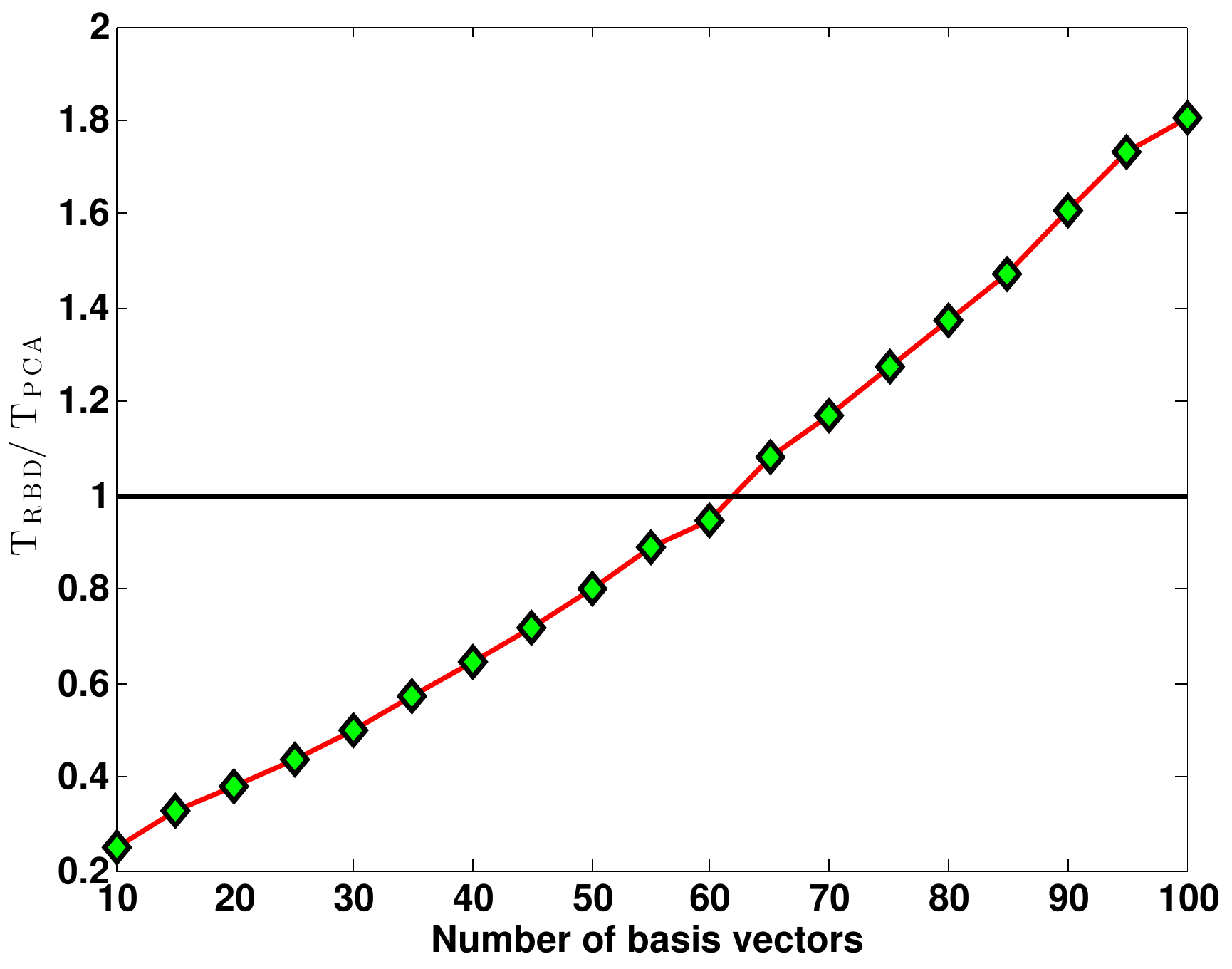}
\caption{Comparison of RBD and other face recognition algorithms. Left: classification error for RBD; Middle: 
classification error for six traditional methods \cite{KokiopoulouChenSaad2011}; Right: Speedup factor of RBD over 
PCA.}
\label{fig:FR_Result}
\end{center}
\end{figure}
We use the UMIST database \cite{GrahamAllinson1998} that is publicly available on Roweis' web page\footnote{http://www.cs.nyu.edu/~roweis/data.html}. 
Table \ref{tab:datasetinfo} summarizes its characteristics:
It contains $20$ people under different poses. The number of different
views per subject varies from $19$ to $48$. We use the cropped version whose snapshot is shown in Figure \ref{fig:UMIST}.

As in \cite{KokiopoulouChenSaad2011}, we randomly choose $10$ views from each class  to form a training set.
The rest of the samples ($375$ of them) are used as testing images. We show the average classification error rates in Figure \ref{fig:FR_Result} Left. 
These averages are computed over $100$ random formations of the training and test sets.
Shown in the middle are the results of six traditional dimension reduction techniques taken from   \cite{KokiopoulouChenSaad2011}. 
Clearly, our method has similar performance as the PCA method, and 
outperforms three of the other five methods. 
However, RBD is much faster than PCA and other methods since they all involves solving eigenproblems \cite{KokiopoulouChenSaad2011}. 
A speedup factor as a function of 
the number of bases is plotted in Figure \ref{fig:FR_Result} Right which 
demonstrates a speedup factor of larger than two for this particular test when we reach the asymptotic region (around 
when the number of basis vectors is $25$).

\section{Concluding remarks}
\label{sec:conclusion}

This paper presents and tests an extremely efficient dimension reduction algorithm for data processing. 
It is multiple times faster than the SVD/PCA-based algorithms. What makes this possible is a greedy algorithm 
that iteratively builds up the reduced space of basis vectors. Each time, the next dimension is located by exploring 
the errors of compression into the current space for all data entries. 
Thanks to an offline-online decomposition mechanism, this searching is independent of the size of each entry. 
Numerical results including one concerning a real world face recognition problem confirm these findings.

\providecommand{\bysame}{\leavevmode\hbox to3em{\hrulefill}\thinspace}
\providecommand{\MR}{\relax\ifhmode\unskip\space\fi MR }
\providecommand{\MRhref}[2]{%
  \href{http://www.ams.org/mathscinet-getitem?mr=#1}{#2}
}
\providecommand{\href}[2]{#2}


\begin{thebibliography}{10}

\bibitem{Almroth_Stern_Brogan}
B.~O. Almroth, P.~Stern, and F.~A. Brogan, \emph{Automatic choice of global
  shape functions in structural analysis}, AIAA Journal \textbf{16} (1978),
  525--528.

\bibitem{Barrault_Nguyen_Maday_Patera}
M.~Barrault, N.~C. Nguyen, Y.~Maday, and A.~T. Patera, \emph{An ``empirical
  interpolation'' method: Application to efficient reduced-basis discretization
  of partial differential equations}, C. R. Acad. Sci. Paris, S\'erie I
  \textbf{339} (2004), 667--672.

\bibitem{BinevCohenDahmenDevorePetrovaWojtaszczyk}
P.~Binev, A.~Cohen, W.~Dahmen, R.~Devore, G.~Petrova, and P.~Wojtaszczyk,
  \emph{Convergence rates for greedy algorithms in reduced basis methods}, SIAM
  J. MATH. ANAL (2011), 1457--1472.

\bibitem{BuffaMadayPateraPrudhommeTurinici2011}
A.~Buffa, Y.~Maday, A.~T. Patera, C.~Prud'homme, and G.~Turinici, \emph{A
  priori convergence of the greedy algorithm for the parametrized reduced
  basis}, ESAIM-Math. Model. Numer. Anal. (2011), Special Issue in honor of
  David Gottlieb.

\bibitem{ChenGottlieb}
Y.~Chen and S.~Gottlieb, \emph{Reduced collocation methods: Reduced basis
  methods in the collocation framework.}, J. Sci. Comput. \textbf{55} (2013),
  no.~3, 718--737.

\bibitem{CHMR_Sisc}
Y.~Chen, J.~S. Hesthaven, Y.~Maday, and J.~Rodr\'{i}guez, \emph{Certified
  reduced basis methods and output bounds for the harmonic {M}axwell's
  equations}, Siam J. Sci. Comput. \textbf{32} (2010), no.~2, 970--996.

\bibitem{CHM_JCP}
Y.~Chen, J.~S. Hesthaven, Y.~Maday, J.~Rodr\'{i}guez, and X.~Zhu,
  \emph{Certified reduced basis method for electromagnetic scattering and radar
  cross section estimation}, CMAME \textbf{233} (2012), 92--108.

\bibitem{Fink_Rheinboldt_1}
J.~P. Fink and W.~C. Rheinboldt, \emph{On the error behavior of the reduced
  basis technique for nonlinear finite element approximations}, Z. Angew. Math.
  Mech. \textbf{63} (1983), no.~1, 21--28. \MR{MR701832 (85e:73047)}

\bibitem{GrahamAllinson1998}
D.~B Graham and N.~M Allinson, \emph{Characterizing virtual eigensignatures for
  general purpose face recognition}, Face Recognition: From Theory to
  Applications (H.~Wechsler, P.~J. Phillips, V.~Bruce, F.~Fogelman-Soulie, and
  T.~S. Huang, eds.), NATO ASI Series F, Computer and Systems Sciences, vol.
  163, 1998, pp.~446--456.

\bibitem{Grepl_Maday_Nguyen_Patera}
M.~A. Grepl, Y.~Maday, N.~C. Nguyen, and A.~T. Patera, \emph{Efficient
  reduced-basis treatment of nonaffine and nonlinear partial differential
  equations}, Mathematical Modelling and Numerical Analysis \textbf{41} (2007),
  no.~3, 575--605.

\bibitem{KokiopoulouChenSaad2011}
E.~Kokiopoulou, J.~Chen, and Y.~Saad, \emph{Trace optimization and
  eigenproblems in dimension reduction methods}, Numerical Linear Algebra with
  Applications \textbf{18} (2011), no.~3, 565--602.

\bibitem{Machielis_Maday_Oliveira_Patera_Rovas}
L.~Machiels, Y.~Maday, I.~B. Oliveira, A.~T. Patera, and D.~V. Rovas,
  \emph{Output bounds for reduced-basis approximations of symmetric positive
  definite eigenvalue problems}, C. R. Acad. Sci. Paris S\'er. I Math.
  \textbf{331} (2000), no.~2, 153--158. \MR{MR1781533 (2001d:65148)}

\bibitem{Maday}
Y.~Maday, \emph{Reduced basis method for the rapid and reliable solution of
  partial differential equations}, International {C}ongress of
  {M}athematicians. {V}ol. {III}, Eur. Math. Soc., Z\"urich, 2006,
  pp.~1255--1270. \MR{2275727 (2007m:65099)}

\bibitem{Maday_Patera_Rovas}
Y.~Maday, A.~T. Patera, and D.~V. Rovas, \emph{A blackbox reduced-basis output
  bound method for noncoercive linear problems}, Nonlinear partial differential
  equations and their applications. Coll\`ege de France Seminar, Vol. XIV
  (Paris, 1997/1998), Stud. Math. Appl., vol.~31, North-Holland, Amsterdam,
  2002, pp.~533--569. \MR{MR1936009 (2003j:65120)}

\bibitem{Maday_Patera_Turinici_2}
Y.~Maday, A.~T. Patera, and G.~Turinici, \emph{A priori convergence theory for
  reduced-basis approximations of single-parameter elliptic partial
  differential equations}, J. Sci. Comput. \textbf{17} (2002), 437--446.

\bibitem{NguyenVeroyPatera2005}
N.C. Nguyen, K.~Veroy, and A.~T. Patera, \emph{Certified real-time solution of
  parametrized partial differential equations}, Handbook of Materials Modeling
  (Sidney Yip, ed.), Springer Netherlands, 2005, pp.~1529--1564 (English).

\bibitem{Noor_Peters}
A.~K. Noor and J.~M. Peters, \emph{Reduced basis technique for nonlinear
  analysis of structures}, AIAA Journal \textbf{18} (1980), no.~4, 455--462.

\bibitem{Porsching}
T.~A. Porsching, \emph{Estimation of the error in the reduced basis method
  solution of nonlinear equations}, Math. Comp. \textbf{45} (1985), no.~172,
  487--496. \MR{MR804937 (86m:65053)}

\bibitem{Prudhomme_Rovas_Veroy_Maday_Patera_Turinici}
C.~Prud'homme, D.~Rovas, K.~Veroy, Y.~Maday, A.~T. Patera, and G.~Turinici,
  \emph{Reliable real-time solution of parametrized partial differential
  equations: Reduced-basis output bound methods}, Journal of Fluids Engineering
  \textbf{124} (2002), no.~1, 70--80.

\bibitem{Rozza_Huynh_Patera}
G.~Rozza, D.B.P. Huynh, and A.T. Patera, \emph{Reduced basis approximation and
  a posteriori error estimation for affinely parametrized elliptic coercive
  partial differential equations: Application to transport and continuum
  mechanics}, Arch Comput Methods Eng \textbf{15} (2008), no.~3, 229--275.

\bibitem{SenNatNorm}
S.~Sen, K.~Veroy, D.B.P. Huynh, S.~Deparis, N.C. Nguyen, and A.T. Patera,
  \emph{``{Natural norm}'' a posteriori error estimators for reduced basis
  approximations}, J. Comput. Phys. \textbf{217} (2006), no.~1, 37 -- 62.

\end{thebibliography}
\end{document}